\documentclass[12pt]{IEEEtran}
\usepackage[utf8]{inputenc}
\usepackage{amssymb}
\usepackage{amsmath}
\usepackage{graphicx}
\usepackage{nomencl}
\usepackage{etoolbox}
\usepackage{algorithm} 
\usepackage{algpseudocode} 
\usepackage{nomencl}
\usepackage[T1]{fontenc}
\usepackage{cite}
\usepackage{array,multirow,graphicx}
\usepackage{multirow}
\usepackage{mathtools}
\usepackage{colortbl}
\usepackage{array}

\title{Flood-aware Optimal Power Flow for Proactive Day-ahead Transmission Substation Hardening}
\author{Mohadese Movahednia, \textit{Student Member, IEEE}, Amin Kargarian, \textit{Senior Member, IEEE}
\date{December 2021}
\thanks{M. Movahednia and A. Kargarian are with the Department of Electrical and Computer Engineering, Louisiana State University, Baton Rouge, LA 70803 USA (e-mail: mmovah4@lsu.edu, kargarian@lsu.edu).}}

\begin{document}
\maketitle
\begin{abstract}
Power system components, particularly electrical substations, may be severely damaged due to flooding, resulting in prolonged power outages and resilience degradation. This problem is more severe in low-elevated regions such as Louisiana. Protective operational actions such as placing tiger dams around substations before flooding can reduce substation vulnerability, damage costs, and energy not supplied cost, thus enhancing power grid resilience. This paper proposes a stochastic mixed-integer programming model for protecting transmission substations one day before flood events using tiger dams. A flood-aware optimal power flow problem is formulated for transmission system operators with respect to the protected/unprotected status of substations. The formulated model aims to identify critical substations whose protection maximizes grid resilience. The output of the proposed model is a day-ahead crew team schedule for installing tiger dams on transmission substations. The effectiveness of the proposed model is evaluated on a 6-bus system, and promising results are obtained.
\end{abstract}
%%%%%%%%%%%%%%%%%%%%%%%%%%%
%%%%%%%%%%%%%%%%%%%%%%%%%%%%%%%%%%%%%%%%%%%%%%%%%%%%%%%%%%%%%
\section{Introduction}
\IEEEPARstart{E}{xtreme} weather events are among the most challenging problems power system operators and utility companies face \cite{panteli2017power, ahmadi2021application}. Long-term and widespread outages, structural damage, human life loss, and multi-million-dollar losses, among others, are unavoidable consequences that are exacerbated by the severity of climate change \cite{bhusal2020power,filabadi2021robust}. In terms of human life loss and component degradation, flooding is the most disastrous climate change for power systems \cite{costa2017substation,AdamsManual1}. Floodwaters cause serious damage to power system infrastructure, notably electrical substations, which are one of the most costly and crucial grid components \cite{karagiannis2017power}.  

A few studies have been conducted recently on improving power system resilience towards flooding. Infrastructure hardening and preemptive action scheduling can mitigate the effects of flooding on the electrical grid \cite{costa2017substation,boggess2014storm,souto2022power}. However, many flooding events have demonstrated that these long-term plans are insufficient, and operational resilience approaches are also essential \cite{wang2015research}. Demand response, energy storage scheduling, and network switching are used in \cite{amirioun2017towards} to reduce load shedding during flooding. To ensure equipment safety, all vulnerable components are tripped out by the operator prior to flooding. Island formation before flooding to minimize load curtailment is presented in \cite{bahrami2021pre}. The impact of floods on a single substation is investigated in \cite{sanchez2020electrical} and \cite{stevens2020interlinking}. The financial cost and feasibility of flood-prevention strategies, such as elevating substations, are assessed in \cite{aerts2014evaluating}. 

However, no operational protective action is suggested in these studies. Protecting electrical substations before flooding by placing inflatable barriers known as tiger dams around them will save millions of dollars and prevent widespread outages. The short installation time of tiger dams, in order of hours, and their low cost have made them practical options. Our previous work presented a stochastic approach for protecting electrical substations by tiger dam one day before the flood \cite{movahednia2021power}. The most critical substations in terms of power grid resiliency and damage costs are identified. Protecting these substations reduces power outages and minimizes power grid cost, including substations' degradation cost and load shedding cost. However, only power distribution substations are considered, and substations are considered independent. Also, network connections and interactions between substations are neglected. These assumptions reduce the effectiveness of the proactive protection scheme.

To fill this gap, this paper presents a decision-making approach to determine optimal day-ahead operational actions to proactively protect transmission substations against flood hazards. A flood-aware optimal power flow (OPF) problem is formulated with respect to substations' protected$/$ unprotected status. The network connection is taken into consideration. Critical substations that might become flooded are identified, and day-ahead crew team scheduling to perform protection actions is obtained. Substations are then protected by installing tiger dams one day before the flood. Simulations on a 6-bus system show that the proposed approach reduces the cost, outage time, and outage magnitude significantly.

%%%%%%%%%%%%%%%%%%%%%%%%%%%%%%%%%%%%%%%%%%%%%%%%%%%%%%%%%%
\section{Proposed Mathematical Model}
A day-ahead stochastic resource scheduling model is formulated to minimize the expected imposed cost to the transmission system due to flooding. Fig. 1 provides an overview of the proposed approach. Mean flood depth, failure probability, repair time, and degradation cost of flood-impacted substations are approximated using historical data and substations' curves, including fragility curve, damage vs. flood depth curve, and repair time vs. damage percentage curve. These data are considered available and provided into the model as input. 

%%%%%%%%%%%%%%%%%%%%%%%%%%%%%%%%%%%%%%%%%%%%%%%%%%%%%%%%%%%%%%%%%%%%%%%%%%%%
\begin{figure}
\centering
    \includegraphics[width=.44\textwidth]{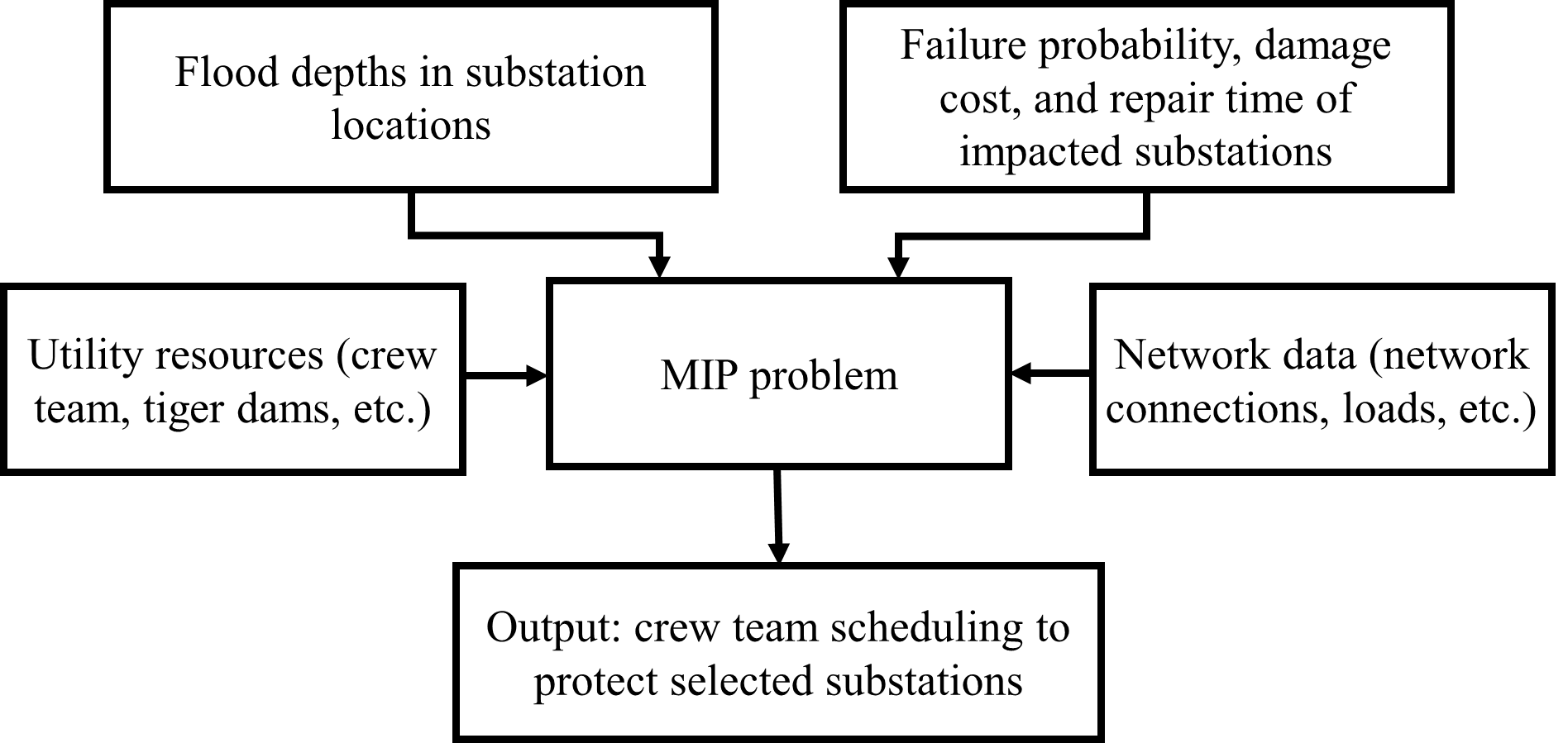}
\caption{An overview of the proposed approach.}
\label{fig1}
\end{figure}
%%%%%%%%%%%%%%%%%%%%%%%%%%%%%%%%%%%%%%%%%%%%%%%%%%%%%%%%%%%%%%%%%%%%%%%%%%%%

\subsection{Objective function}
The objective function is formulated in (1). The proposed model minimizes the expected cost imposed on the transmission system due to substation failure and energy not supplied.

\begin{equation} \tag{1}
\label{1}
\min_{} \sum_{\forall k}\sum_{\forall s} \pi_s F_{ks}  \beta_k + \sum_{\forall s} \pi_s \sum_{\forall i} \sum_{\forall t} VOLL\times \mathcal {LS}_{ist}
\end{equation}
where $VOLL$ is the value of lost load. $\mathcal {LS}_{ist}$ is load shedding of bus $i$ at time $t$ for scenario $s$. $\pi_s$ is the probability of scenario $s$. $F_{ks}$ is substation failure scenarios that indicates whether substation $k$ is failed or not in scenario $s$. $F_{ks}=1$ means that substation $k$ will fail due to flooding, and zero means it will survive. The probability of scenario $s$ is calculated as follows \cite{movahednia2021power}:

\begin{equation} \tag{2}
\label{2}
\pi_s=\prod_{\forall k} \left(F_{ks} \pi_k +\left(1-F_{ks}\right) \left(1-\pi_k \right) \right);  \mbox{  } \forall s.
\end{equation}
where $\pi_k$  is the failure probability of substation $k$. Scenario probabilities are calculated as the product of the failure probability of failed substations and the survival probability of survived substations. The failure probability of failed substation $k$ is $F_{ks} \pi_k$. The value of $F_{ks}$ for a failed substation is one. Thus, we have $F_{ks} \pi_k=1\times\pi_k=\pi_k $. The survival probability of survived substation $k$ is $\left(1-F_{ks}\right) \left(1-\pi_k \right)$, which is $\left(1-0\right)\times\left(1-\pi_k \right)=1-\pi_k$ since $F_{ks}=0$ for a survived substation.
\subsection{Constraints}
In (1), $\beta_k$ is the cost imposed on transmission substation $k$ due to the flood. 

\begin{equation} \tag{3}
\label{3}
\beta_k=\theta_k \mathcal{C}_k+\left(1-\theta_k\right)\mathcal{DC}_k;  \mbox{  } \forall k. 
\end{equation}
where $\theta_k$ is a binary variable that shows whether substation $k$ is selected for protection or not. $\mathcal{C}_k$ denotes the cost of tiger dams for protecting substation $k$. $\mathcal{DC}_k$ is the cost of structural damages to substation $k$. $\mathcal{DC}_k$ is estimated according to flood depth estimation, damage percentage, and substation price.

{\normalfont\itshape Network, Substation, and Generator Constraints:} The availability of substation $k$ in scenario $s$, denoted by $\mathfrak{h}_{ks}$, depends on its failure and protection status and is calculated by (4). 
\begin{equation} \tag{4}
\label{4}
\mathfrak{h}_{ks}=1-F_{ks}\left(1-\theta_k\right);\mbox{  } \forall k,\forall s.
\end{equation}

The power balance equation for each bus $i$ is:
\begin{equation} \tag{5}
\label{5}
\sum_{\forall g} p_{gst}+\mathcal {LS}_{ist}-P_{ist}^d=\sum_{\forall l}f_{lst};\mbox{  } \forall i,\forall s,\forall t.
\end{equation}
where $p_{gst}$, $P_{ist}^d$, and $f_{lst}$ are, respectively, power generation of generator $g$ located at node $i$, load connected to node $i$, and power flow through lines connected to node $i$ for scenario $s$ at time $t$. Load shedding is limited by (6).
\begin{equation} \tag{6}
\label{6}
0\leq \mathcal {LS}_{ist}\leq P_{ist}^d;\mbox{  } \forall i,\forall s,\forall t.
\end{equation}

The power flow in transmission line $l$ is modeled by (7)-(10). The availability of transmission lines depends on the availability of their end substation terminals. If both substation terminals of line $l$ are available, $\mathfrak{h}_{ls}$ equals to 1 and $f_{lst}=B_l \left(\triangle \delta_{st}\right)$. $\mathfrak{h}_{ls}=0$ if at least one substation terminal is out; thus $f_{lst}$ would be zero and independent of voltage phase angles of terminal buses. 
\begin{equation} \tag{7}
\label{7}
-\mathfrak{h}_{ls} f_{l}^{max}\leq f_{lst}\leq \mathfrak{h}_{ls} f_{l}^{max};\mbox{  } \forall l,\forall s,\forall t.
\end{equation}
\begin{equation} \tag{8}
\label{8}
f_{lst}+ M \left(1-\mathfrak{h}_{ls}\right) \geq B_l \left(\triangle \delta_{st}\right);\mbox{  } \forall l,\forall s,\forall t.
\end{equation}
\begin{equation} \tag{9}
\label{9}
f_{lst} \leq M \left(1-\mathfrak{h}_{ls}\right) + B_l \left(\triangle \delta_{st}\right);\mbox{  } \forall l,\forall s,\forall t.
\end{equation}
\begin{equation} \tag{10}
\label{10}
\mathfrak{h}_{ls}=\mathfrak{h}_{k_o s} \mathfrak{h}_{k_d s};\mbox{  } \forall l,\forall s.
\end{equation}
where $f_{l}^{max}$ is the capacity of line $l$. $M$ is a big  number. $B_l$ is susceptance of line $l$. $\triangle \delta_{st}$ is the voltage angle difference between the origin and destination nodes of line $l$ in scenario $s$. $\mathfrak{h}_{k_o s}$ and $\mathfrak{h}_{k_d s}$ are availability status of origin and destination substations of line $l$.

The upper and lower bounds of thermal generation units, which depend on the availability of substation $k$, are enforced by (11). The ramp up and ramp down constraints of thermal units are presented in (12) and (13).

\begin{equation} \tag{11}
\label{11}
\mathfrak{h}_{ks} P_g^{min} \leq p_{gst} \leq \mathfrak{h}_{ks} P_g^{max};\mbox{  } \forall g,\forall s,\forall t.
\end{equation}
\begin{equation} \tag{12}
\label{12}
p_{gst} - p_{gst-1} \leq \mathcal{RU}_{g};\mbox{  } \forall g,\forall s,\forall t.
\end{equation}
\begin{equation} \tag{13}
\label{13}
p_{gst-1} - p_{gst} \leq \mathcal{RD}_{g};\mbox{  } \forall g,\forall s,\forall t.
\end{equation}

{\normalfont\itshape Crew Team Constraints:} The sufficient time to place tiger dams around substations, $\tau_k$, is approximated by (14). The flood depth and available crew members determine the required time. 
\begin{equation} \tag{14}
\label{14}
\tau_k = {4+10 \left(\mu_k -0.45 \right)\over \mathcal{N}_{c}^{team}};\mbox{  } \forall 0.45 \leq  \mu_k \leq 1.5.
\end{equation}
where $\mu_k$ is the mean value of flood depth, and $\mathcal{N}_{c}^{team}$ is the number of members per team.

The scheduled program of crew teams is constrained by (15)-(20). Index $\tilde{t}$ shows the time intervals of protection actions. Binary variable $x_{nk\tilde{t}}$ shows the status of crew team $n$ at substation $k$ at time $\tilde{t}$. If $x_{nk\tilde{t}}$ is one, crew team $n$ is scheduled to mount tiger dam around substation $k$ at time $\tilde{t}$. Equation (15) enforces the limitation of total available time and crew teams. Equation (16) ensures that the total time spent by crew teams at substation $k$ is equal to $\tau_k$. Inequality (17) guarantees that each crew team can only perform tiger dam installation at one substation at a time. One crew team must install tiger dams around substation $k$ in sequential hours. Constraints (18) and (19) prevent multiple crew dispatch to one substation. Binary variable $y_{nk\tilde{t}}$ becomes one if crew team $n$ is sent to substation $k$ at time $\tilde{t}$. Equation (20) limits the summation of $y_{nk\tilde{t}}$ over time and crew teams to prevent multiple teams dispatch to substation $k$. 
\begin{equation} \tag{15}
\label{15}
\sum_{\forall n} \sum_{\forall k} \sum_{\forall \tilde{t}} x_{nk \tilde{t}}  \leq \mathcal{N}^{team} \mathcal{T}
\end{equation}
\begin{equation} \tag{16}
\label{16}
\sum_{\forall n} \sum_{\forall \tilde{t}} x_{nk \tilde{t}} = \theta_k \tau_k ;\mbox{  } \forall k.
\end{equation}
\begin{equation} \tag{17}
\label{17}
\sum_{\forall k} x_{nk \tilde{t}} \leq 1;\mbox{  } \forall n, \forall \tilde{t}.
\end{equation}
\begin{equation} \tag{18}
\label{18}
{x_{nk \tilde{t}}-x_{nk \tilde{t}-1} \over 2}- \epsilon  \leq y_{nk\tilde{t}};\mbox{  }  \forall n, \forall k, \forall \tilde{t}.
\end{equation}
\begin{equation} \tag{19}
\label{19}
y_{nk\tilde{t}} \leq {1+x_{nk \tilde{t}}-x_{nk \tilde{t}-1} \over 2}+ \epsilon;\mbox{  }  \forall n, \forall k, \forall \tilde{t}.
\end{equation}
\begin{equation} \tag{20}
\label{20}
\sum_{\forall n} \sum_{\forall \tilde{t}} y_{nk\tilde{t}} \leq 1;\mbox{  }  \forall k.
\end{equation}

\subsection{Constraint Linearization}
In (10), $\mathfrak{h}_{ls}$ is the product of two binary variables and is linearized as follows:
\begin{equation} \tag{21}
\label{21}
\mathfrak{h}_{ls} \leq \mathfrak{h}_{k_o s};\mbox{  }  \forall l, \forall s .
\end{equation}
\begin{equation} \tag{22}
\label{22}
\mathfrak{h}_{ls} \leq \mathfrak{h}_{k_d s};\mbox{  }  \forall l, \forall s .
\end{equation}
\begin{equation} \tag{23}
\label{23}
\mathfrak{h}_{ls} \geq \mathfrak{h}_{k_o s} + \mathfrak{h}_{k_d s} - 1;\mbox{  }  \forall l, \forall s .
\end{equation}

Standard MIP solvers can solve the resultant optimization model.  
\subsection{Model Summary}
The proposed model is a mixed-integer linear program as follows:
\begin{equation} \tag{24}
\label{24}
\min_{}  (1)
\end{equation}
\centerline{s.t. (\ref{2})-(\ref{9}), (\ref{11})-(\ref{23})}

%%%%%%%%%%%%%%%%%%%%%%%%%%%%%%%%%%%%%%%%%%%%%%%%%%%%%%%%%%%%%
\section{Numerical Simulation}
The effectiveness of the proposed model is evaluated on a 6-bus system shown in Fig. 2. Bus data are given in Table I. Substation information is given in Table II. It is assumed that substation data, including flood depth, failure probability, expected degradation cost, and expected repair time due to flooding, is approximated using historical data and is fed to the proposed model as input. Two teams, each with four members, are available. The total time available to perform protective actions a day before flooding is considered five hours. Four representative substation failure scenarios are generated using the presented scenario reduction method in [16]. Substation failure scenarios and their probabilities are presented in Table III.

The proposed approach seeks to reduce the expected damage cost and energy not supplied cost of the whole system. In this regard, flood-prone substations with a higher failure likelihood and cost are more likely to be protected. The crew team schedules are presented in Tables IV. The four substations selected for protection are $k2$, $k3$, $k5$, and $k6$. Substation $k6$ has the highest failure probability. Also, substation $k2$ has a high failure probability, and protecting it would reduce the expected cost. The substations $k5$ and $k3$ have relatively high failure rates. Moreover, they are load terminals, and protecting them would reduce the expected load shedding. Substation $k4$ has a high failure probability, but it is not selected for protection. As shown in Fig. 2, there are two generators in the system. According to Table I, the generator connected to substation $k1$, which has a low failure probability, can supply the total load of the system. In this case, losing the connected generator to substation $k4$ does not result in significant load shedding in the system. Therefore, protecting substation $k4$ is given a lower priority than other substations.

The results of the proposed stochastic approach and a deterministic substation protection approach are compared. In the deterministic approach, the failure probability of substations is ignored, and the substation protection decisions are made assuming that all substations will fail due to flood. Based on the considered four substation failure scenarios, expected cost, power outage magnitude and time are determined to compare the deterministic and stochastic approaches fairly. The resilience curves of these methods are illustrated in Fig. 3. According to Fig. 3, the expected power outage magnitude and outage time of the stochastic approach are, respectively, 67\% and 11\% less than the deterministic approach. Protecting more critical substations by the stochastic approach, more resilient power grid is achieved as less power outage is experienced, and the system returns to a normal status faster. Table V shows the expected cost imposed on the transmission system due to flood, which is 53\% less in the case of the proposed stochastic approach. 

%%%%%%%%%%%%%%%%%%%%%%%%%%%%%%%%%%%%%%%%%%%%%%%%%%%%%%%%%%%%%%%%%%%%%%%%%%%%
\begin{figure}
\centering
    \includegraphics[width=.4\textwidth]{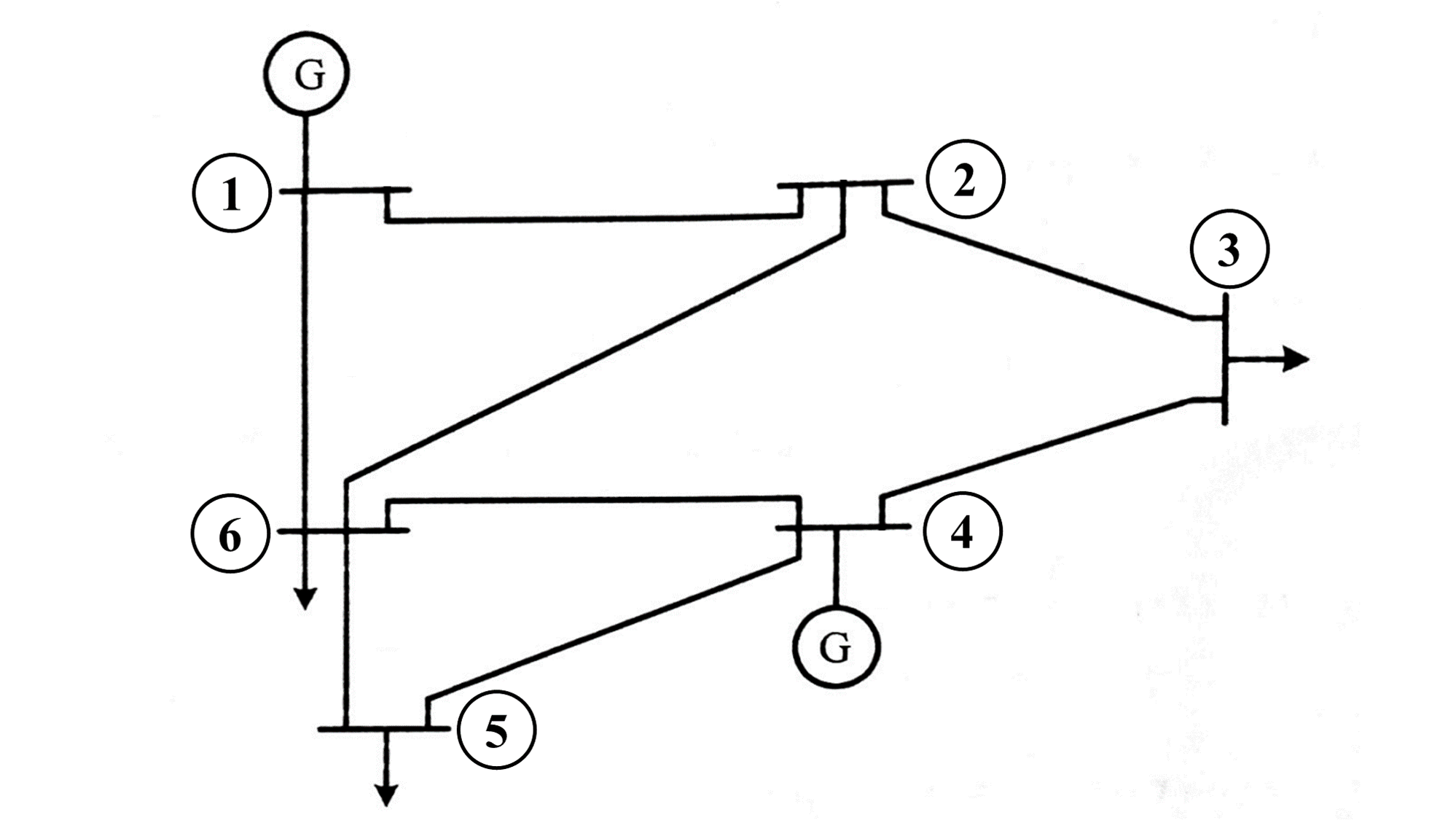}
\caption{6-bus system}
\label{fig2}
\end{figure}
%%%%%%%%%%%%%%%%%%%%%%%%%%%%%%%%%%%%%%%%%%%%%%%%%%%%%%%%%%%%%%%%%%%%%%%%%%%%

%%%%%%%%%%%%%%%%%%%%%%%%%%%%%%%%%%%%%%%%%%%%%%%%%%%%%%%%%%
%%%%%%%%%%%%%%%%%%%%%%%%%%%%%%%%%%%%%%%%%%%%%%%%%%%%%%%%%%%
\begin{table}
\setlength{\arrayrulewidth}{0.2mm}
\begin{center}
\caption{6-Bus System Data}
\label{table:1}
\begin{tabular}{c  c c c c c } \hline 
Bus & Average demand & $P_g^{min}$ & $P_g^{max}$ &$\mathcal{RU}_{g}$ &$\mathcal{RD}_{g}$\\ \hline\hline
1  &- &50    &150  &14   &14\\ \hline
3  &55  &-  &- &-   &- \\ \hline
4 &-  &30  & 100 & 10 &10  \\ \hline
5 &30  &- &- &- &- \\ \hline
6 &50  &- &-  &- &- \\
 \hline
\end{tabular}
\end{center}
\end{table}
%%%%%%%%%%%%%%%%%%%%%%%%%%%%%%%%%%%%%%%%%%%%%%%%%%%%%%%%%%
%%%%%%%%%%%%%%%%%%%%%%%%%%%%%%%%%%%%%%%%%%%%%%%%%%%%%%%%%%%%%%%%%%%%%%%%%%%%

%%%%%%%%%%%%%%%%%%%%%%%%%%%%%%%%%%%%%%%%%%%%%%%%%%%%%%%%%%
%%%%%%%%%%%%%%%%%%%%%%%%%%%%%%%%%%%%%%%%%%%%%%%%%%%%%%%%%%%

\begin{table*}
\setlength{\arrayrulewidth}{0.2mm}
\begin{center}
\caption{Substation Information}
\label{table:2}
\begin{tabular}{c c c c c c } \hline 
Sub. & Mean flood depth (m) & Failure probability & Repair time (h) &Damage cost (\$) &Protection time (h)\\ \hline\hline
 k1  &0.5  &0.115  &12.3  &40410   &2\\ \hline
 k2  &0.9  &0.276  &19.4  &65464   &3\\ \hline
 k3  &0.8  &0.235  &17.7  &59403   &2\\ \hline
 k4  &1.1  &0.361  &23.1  &77809   &3 \\ \hline
 k5  &0.6  &0.154  &14.1  &46890   &2\\ \hline
 k6  &1.2  &0.404  &25    &83960   &3\\ 
 \hline
\end{tabular}
\end{center}
\end{table*}

%%%%%%%%%%%%%%%%%%%%%%%%%%%%%%%%%%%%%%%%%%%%%%%%%%%%%%%%%%
%%%%%%%%%%%%%%%%%%%%%%%%%%%%%%%%%%%%%%%%%%%%%%%%%%%%%%%%%%%%%%%%%%%%%%%%%%%%

%%%%%%%%%%%%%%%%%%%%%%%%%%%%%%%%%%%%%%%%%%%%%%%%%%%%%%%%%%
%%%%%%%%%%%%%%%%%%%%%%%%%%%%%%%%%%%%%%%%%%%%%%%%%%%%%%%%%%%
\begin{table}
\setlength{\arrayrulewidth}{0.2mm}
\begin{center}
\caption{Substation Failure Scenarios}
\label{table:3}
\begin{tabular}{m{1.1cm}  c c c c c c c } \hline 
Sce./Sub. & k1 & k2 & k3 &k4 &k5 &k6 &Probability \\ \hline\hline
 S1  & &    &  &   &  &Fail &0.612\\ \hline
 S2  &  &  & & Fail  & &Fail &0.346\\ \hline
 S3  &  & Fail  & Fail & Fail &  &Fail &0.041 \\ \hline
 S4  & Fail & Fail  & Fail & Fail & Fail &Fail &9.6e-4 \\
 \hline
\end{tabular}
\end{center}
\end{table}
%%%%%%%%%%%%%%%%%%%%%%%%%%%%%%%%%%%%%%%%%%%%%%%%%%%%%%%%%%
%%%%%%%%%%%%%%%%%%%%%%%%%%%%%%%%%%%%%%%%%%%%%%%%%%%%%%%%%%%%%%%%%%%%%%%%%%%%
%%%%%%%%%%%%%%%%%%%%%%%%%%%%%%%%%%%%%%%%%%%%%%%%%%%%%%%%%%
%%%%%%%%%%%%%%%%%%%%%%%%%%%%%%%%%%%%%%%%%%%%%%%%%%%%%%%%%%%
\begin{table}
\setlength{\arrayrulewidth}{0.2mm}
\begin{center}
\caption{Crew Team Schedule}
\label{table:4}
\begin{tabular}{m{1.8cm}  m{0.9cm} m{0.9cm} m{0.9cm} m{0.9cm} m{0.9cm}} \hline 
Team/time (h) &1 &2 &3 &4 &5\\ \hline\hline
\multirow{2}{*}{Team 1}       &\multicolumn{2}{c}{\cellcolor{yellow}k5}                                             &\multicolumn{3}{c}{}\\
                  &\multicolumn{2}{c}{}   &\multicolumn{3}{c}{\cellcolor{yellow}k2} \\ \hline
\multirow{2}{*}{Team 2}       & \multicolumn{3}{c}{\cellcolor{green} k6}                                              &\multicolumn{2}{c}{}\\ 
                  &\multicolumn{3}{c}{}   &\multicolumn{2}{c}{\cellcolor{green}k3} \\ \hline
\end{tabular}
\end{center}
\end{table}
%%%%%%%%%%%%%%%%%%%%%%%%%%%%%%%%%%%%%%%%%%%%%%%%%%%%%%%%%%
%%%%%%%%%%%%%%%%%%%%%%%%%%%%%%%%%%%%%%%%%%%%%%%%%%%%%%%%%%%%%%%%%%%%%%%%%%%%
%%%%%%%%%%%%%%%%%%%%%%%%%%%%%%%%%%%%%%%%%%%%%%%%%%%%%%%%%%%%%%%%%%%%%%%%%%%%
\begin{figure}
\centering
    \includegraphics[width=.42\textwidth]{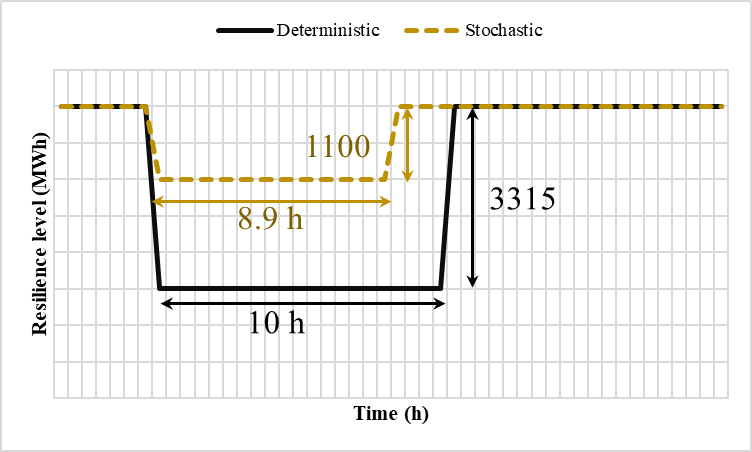}
\caption{Resilience curve of stochastic and deterministic approaches.}
\label{fig3}
\end{figure}
%%%%%%%%%%%%%%%%%%%%%%%%%%%%%%%%%%%%%%%%%%%%%%%%%%%%%%%%%%%%%%%%%%%%%%%%%%%%
%%%%%%%%%%%%%%%%%%%%%%%%%%%%%%%%%%%%%%%%%%%%%%%%%%%%%%%%%%
%%%%%%%%%%%%%%%%%%%%%%%%%%%%%%%%%%%%%%%%%%%%%%%%%%%%%%%%%%%
\begin{table}[h!]
\setlength{\arrayrulewidth}{0.2mm}
\begin{center}
\caption{Expected Cost Of Deterministic And Stochastic Models}
\label{table:5}
\begin{tabular}{m{3.6cm}  m{3.6cm} } \hline 
Model &Expected cost (\$)\\ \hline\hline
Proposed stochastic  &85,164 \\ \hline
Deterministic  &181,119  \\ \hline
\end{tabular}
\end{center}
\end{table}
%%%%%%%%%%%%%%%%%%%%%%%%%%%%%%%%%%%%%%%%%%%%%%%%%%%%%%%%%%
%%%%%%%%%%%%%%%%%%%%%%%%%%%%%%%%%%%%%%%%%%%%%%%%%%%%%%%%%%%%%%%%%%%%%%%%%%%%
%%%%%%%%%%%%%%%%%%%%%%%%%%%%%%%%%%%%%%%%%%%%%%%%%%%%%%%%%%
%%%%%%%%%%%%%%%%%%%%%%%%%%%%%%%%%%%%%%%%%%%%%%%%%%%%%%%%%%%%%%%%%%%%%%%%%%%%
\section{Conclusion}
Flooding imposes severe stress on the power system, especially electrical substations. The extensive and prolonged outages imposed by flooding can be alleviated by performing operational protective actions, such as protecting electrical substations by tiger dams. In this regard, determining critical substations and allocating available resources are crucial. This paper presents a mixed-integer linear optimization model that minimizes the cost imposed on the transmission system due to flooding. A stochastic flood-aware OPF is formulated. The proposed approach identifies critical transmission substations that should be protected one day before flooding and determines the scheduling program of crew teams. Simulation results on the 6-bus systems show that the expected cost, power outage magnitude, and outage duration are decreased respectively by 53\%, 67\%, and 11\% as compared to a deterministic approach.
%%%%%%%%%%%%%%%%%%%%%%%%%%%%%%%%%%%%%%%%%%%%%%%%%%%%%%%%%%%%%
\bibliographystyle{ieeetr}
\bibliography{FloodawareOPF.bib}
%%%%%%%%%%%%%%%%%%%%%%%%%%%%%%%%%%%%%%%%%%%%%%%%%%%%%%%%%%%%%
\end{document}